%Paper: alg-geom/9502010
%From: Gaitsgory Denis <gaitsgde@math.tau.ac.il>
%Date: Tue, 14 Feb 1995 15:28:28 +0200
%Date (revised): Wed, 15 Feb 1995 14:09:41 +0200

\documentstyle{amsppt}

%\font\ref=cmr9
%\font\refit=cmti9
%\font\refbf=cmbx9

\nopagenumbers

\input amstex
\magnification=\magstep1

%Special fonts
\font\tenboldgreek=cmmib10  \font\sevenboldgreek=cmmib10 at
7pt
\font\fiveboldgreek=cmmib10 at 7pt
\newfam\bgfam
\textfont\bgfam=\tenboldgreek \scriptfont\bgfam=\sevenboldgreek
\scriptscriptfont\bgfam=\fiveboldgreek

%\def\bg{\fam6}
%\mathchardef\alpha="700B
%\def\bfalp{{\fam=\bgfam\balp}}

\font\tengerman=eufm10 \font\sevengerman=eufm7 \font\fivegerman=eufm5
\font\tendouble=msym10 \font\sevendouble=msym7 \font\fivedouble=msym5

\textfont4=\tengerman \scriptfont4=\sevengerman
\scriptscriptfont4=\fivegerman
\newfam\dbfam
\textfont\dbfam=\tendouble \scriptfont\dbfam=\sevendouble
\scriptscriptfont\dbfam=\fivedouble
\def\gr{\fam4}

\mathchardef\ng="702D
\mathchardef\dbA="7041
\mathchardef\sm="7072
\mathchardef\nvdash="7030
\mathchardef\nldash="7031
\mathchardef\lne="7008
\mathchardef\sneq="7024
\mathchardef\spneq="7025
\mathchardef\sne="7028
\mathchardef\spne="7029
\mathchardef\ltms="706E
\mathchardef\tmsl="706F

\mathchardef\dbA="7041

\def\sdp{\times \hskip -0.3em {\raise 0.3ex
\hbox{$\scriptscriptstyle |$}}} % semidirect product

\def\subheading#1{\medskip\goodbreak\noindent{\bf
#1.}\quad}

\def\longmapright #1 #2 {\smash{\mathop{\hbox to
#1pt {\rightarrowfill}}\limits^{#2}}}
\def\longmapleft #1 #2 {\smash{\mathop{\hbox to
#1pt {\leftarrowfill}}\limits^{#2}}}

\def\back{{\raise 2.5pt\hbox{$\,\scriptscriptstyle\backslash\,$}}}

\def\part{\partial}
\def\lwr #1{\lower 5pt\hbox{$#1$}\hskip -3pt}
\def\rse #1{\hskip -3pt\raise 5pt\hbox{$#1$}}
\def\lwrs #1{\lower 4pt\hbox{$\scriptstyle #1$}\hskip -2pt}
\def\rses #1{\hskip -2pt\raise 3pt\hbox{$\scriptstyle #1$}}

\def\<#1{\left\langle{#1}\right\rangle}

\def\subinbn{{\subset\hskip-8pt\raise
0.95pt\hbox{$\scriptscriptstyle\subset$}}}

\def\llvdash{\mathop{\|\hskip-2pt \raise 3pt\hbox{\vrule
height 0.25pt width 1.5cm}}}

\def\lvdash{\mathop{|\hskip-2pt \raise 3pt\hbox{\vrule
height 0.25pt width 1.5cm}}}

\def\fakebold#1{\leavevmode\setbox0=\hbox{#1}%
  \kern-.025em\copy0 \kern-\wd0
  \kern .025em\copy0 \kern-\wd0
  \kern-.025em\raise.0333em\box0 }

\font\msxmten=msxm10
\font\msxmseven=msxm7
\font\msxmfive=msxm5
\newfam\myfam
\textfont\myfam=\msxmten

\scriptfont\myfam=\msxmseven
\scriptscriptfont\myfam=\msxmfive
\mathchardef\rhookupone="7016
\mathchardef\ldh="700D
\mathchardef\leg="7053
\mathchardef\ANG="705E
\mathchardef\lcu="7070
\mathchardef\rcu="7071
\mathchardef\leseq="7035
\mathchardef\qeeg="703D
\mathchardef\qeel="7036
\mathchardef\blackbox="7004
\mathchardef\bbx="7003
\mathchardef\simsucc="7025

\def\+{\oplus }

\define\trrt{\Upsilon(T,t')\underset{t\neq t'}\to\bigotimes
O(S_t)\otimes\Upsilon(S_{t'},s')}
\redefine\Lamda{\wedge}
\define\arss{A^{\otimes S-s'}}
\define\Ost{\underset{t\in T}\to\bigotimes O(S_t)}
\define\Rm{R\text{-modules }}
\define\sre{\sum_{s\in S}}
\redefine\pm{{\rm pseudo-monoidal  }}
\define\Ei{{Ext_{A\otimes A}}^{i+1}}
\define\for{\text{for each finite set}\ S\text{ }}
\define\sbh{\subheading}
\define\fs{F^C_S}
\define\CSo{(C^o)^{\times S}}
\define\st{S\to T\text{ is a surjection of finite sets}}
\redefine\bt{\underset{t\in T}\to\times B_t}
\define\CT{C^{\times T}}
\define\ft{F^C_T}
\define\fst{F^C_{S_t}}
\define\fsst{F^{\Upsilon}_{S_{t'},s'}}
%\define\fors{\text{for any finite set } S\text{\and } s\in S}
\define\fss{F^{\Upsilon}_{S,s'}}
\define\Csso{(C^o)^{\times (S-s')}}
\define\sst{S=\underset{t\in T}\to\bigcup S_t,s'\in S_{t'}}
\define\rti{R[t]/t^{i+1}\cdot R[t]}
\define\Cai{Deform^{i+1}_{A_i}(A)}
\define\g{{\gr g}}
\define\fm{{\Im}_M}
\define\fa{{\Im}_A}
\define\as{A^{\otimes S}}
\redefine\Tau{\Upsilon}
\define\pbw{Poincar\'e-Birkhoff-Witt theorem }
\redefine\RR{{\Cal R}}
\redefine\LL{{\Cal L}}
\redefine\FF{{\Cal F}}
\redefine\tau{\upsilon}
\redefine\un{1}

\topmatter
\title Operads, Grothendieck topologies \\
and deformation theory
\endtitle
\author  Dennis Gaitsgory \endauthor
\address{School of Mathematical Sciences, Tel-Aviv University,
Ramat-Aviv, Israel} \endaddress
\email
gaitsgde\@math.tau.ac.il
\endemail
\endtopmatter

\heading 0. Introduction \endheading
\sbh{0.1}
Gerstenhaber's papers in the Annals showed that deformation theory of
associative algebras over a field is "controlled" by Hochschild cohomology.
The passage from deformations to cohomology is realised by
means of Hochschild cochains. This approach has two main drawbacks:
\sbh{1} It is impossible to generalize it to the case of
algebras over a ring, or more generally, to sheaves of algebras
over a scheme. This is because neither deformations are described
by cochains, nor cohomologies can be computed using bar-resolution.
\sbh{2} Although deformations and cohomology are
invariant objects they are connected by choosing some specific resolution.
\sbh{0.2}
Our aim in the present work is to define an appropriate cohomology
theory and to find an invariant way to pass from deformations to
cohomology.

The initial idea, which goes back, probably, to Quillen, is to describe
deformations of an algebraic object (e.g. associative algebra) by means of
"resolving" it by free objects of the same type (in our example, free
associative algebras). In principal, all the results can be formulated
already on the level of our initial algebra $A$ and a free algebra $B$
mapping onto it. However, the picture is much easier to grasp,
when we consider the category of all algebras over $A$. These algebras
form a site, and cohomologies that we are looking for are just cohomologies
of certain sheaves on this site. This is the main idea of the paper.

An advantage of this approach is that we can treat in the same framework
algebras of all types, i.e. algebras over an arbitrary operad.
\sbh{0.3}
Let us describe briefly the contents of the paper.

In sections 1 and 2 we describe the formalism of operads, algebras
over operads and modules over them. Our presentation is inspired
by some ideas of A.Beilinson [6] and is very close to that of [1].
The essential difference is that we are using the language of
\pm categories.

In section 3 we introduce the site $C(A)$ and study the connection
between the category of sheaves on this site and the category
of $A$-modules. In paricular, we introduce the notion of the cotangent
complex of an algebra. The site $C(A)$ was introduced
first in [1]. However, one can think of deformation theory
(e.g. Theorem 4.2) as giving a hint how to define correct cohomologies:
just look at deformations of the corresponding object.

In section 4 we study deformations of an algebra over
an operad. Cohomology classes that arise in deformation theory
are realised as classes of certain torsors and gerbes.

Finally, in section 5 we give an application of the theory
presented in the paper. We prove the \pbw for Lie algebras
which are flat modules over the ground ring. The idea of such an
approach to the \pbw belongs to J.Bernstein and has
been already realized in [6].
\sbh{Acknowledgements}
The author would like to thank A.Joseph, J.Bernstein, V.Hinich, S.Shnider,
L.Breen, A.Braverman and R.Bezrukavnikov for interesting and
stimulating discussions.

\heading 1. Pseudo-monoidal categories and operads  \endheading
\subheading{1.0}
We are working over a fixed ground ring $R$ and all categories are assumed
to be $R$-linear. If $C$ is a gategory, $C^o$ will denote the
opposite category.
\subheading{1.1 Pseudo-monoidal categories}
Let $C$ be a category. A pseudo-monoidal structure on it is a collecion
of functors $\for$:

$\fs:\CSo\times C\to \Rm$ equipped with the
following additional data (composition maps):

If $\st$, for each element $\bt$ of $\CT$ we are given a natural
transformation
between two functors $\CSo\times C\to \Rm$:
$${\bigotimes}_{t\in T}\fst(\underset{s\in S_t}\to\times A_s,B_t)
\otimes\ft(\underset{t\in T}\to\times
B_t,D)
\longrightarrow \fs(\underset{s\in S}\to\times A_s,D),$$
($S_t$ denotes here the preimage of $t\in T$)
with $F_1^C(A,B)=Hom(A,B)$ (subscript $1$ means one element set)
such that these natural transformations are compatible with respect to
compositions of partitions in the obvious sense.
\subheading{1.1.1 Example}
If $C$ is a strictly symmetric monoidal category, we endow it with a
pseudo-monoidal structure by setting
$$\fs(\underset{s\in S}\to\times A_s,B)=Hom(\otimes_{s\in S}A_s,B)$$
where the above natural transformations are just composition maps.
\sbh{1.1.2 Variant}
As in usual monoidal categories, one can require existence of a unit
object in a \pm category. This means that there must be an object
$\un\in C$ such that $\for$, for each $\underset{s\in S}\to\times
A_s\in C^{\times S}$ and for each $B\in C$, $F^C_{S\cup 1}(1\underset{s\in S}
\to\times A_s,B)$ is canonically isomorphic to $F^C_S(\underset{s\in S}
\to\times A_s,B)$. Such \pm categories will be called unital.

\subheading{1.1.3}
A pseudo-monoidal functor between two pseudo-monoidal
categories $C_1$ and $C_2$ is a (covariant) functor $G: C_1\to C_2$
equipped with a
natural transformation $F^{C_1}_S\to F^{C_2}_S\circ G$ which
is compatible with composition maps $\for$.
Pseudo-monoidal natural transformations are defined in a similar way.

\subheading{1.2 Operads}
Operad is by definition a \pm category $O$ which has essentially
one object.
\subheading{1.2.1}
Equivalently,
one can view  operads as the following linear algebra data:
\roster
\item A collection of $\Rm$ \for denoted by $O(S)$ (thought of
as $\fs(A^{\times S},A)$ for $A\in C$)
\item A distinguished element $1\in O(1)$
\item If $\st$ there is a composition map \newline
$$O(T)\underset{t\in T}\to\bigotimes O(S_t)\to O(S)$$
\endroster
that satisfies
\roster
\item The composition $1\otimes O(S)\to O(1)\otimes O(S)\to O(S)$ \newline
is the identity map  $\for$.
\item Composition maps are compatible with compositions of partitions.
\endroster
\sbh{1.2.2 Variant}
We define unital operads as unital \pm categories having at most
one isomorphism class of objects distinct from $\un$. It is
not difficult to work out this definition also in linear
algebra terms.

\sbh{1.2.3}
If $R\to R'$ is a ring homomorphism, to any operad over $R$
one can assign an operad over $R'$ by taking
tensor products over $R$ with $R'$. We will denote them
by same letters when no confusion can occur.

\subheading{1.3 Algebras over an operad}
In order to simplify the exposition we will consider
only algebras of $\Rm$.
In principal, one can define algebras over an operad in any strictly
symmetric monoidal category and develop deformation theory for them.

An $O$-algebra (of $\Rm$), or equivalently, an algebra over $O$
is by definition a \pm functor $O\to \Rm$.
Morphisms between $O$-algebras are defined
to be \pm natural transformations between such functors.
\subheading{1.3.1}
In linear agebra terms, an $O$-algebra is
an $R$-module $A$ which \newline
$\for$ is endowed with a map $$O(S)\otimes A^{\otimes S}\to A$$
such that if $\st$ the square
$$ \CD
O(T)\Ost\otimes A^{\otimes S}   @>>>    O(S)\otimes A^{\otimes S}  \\
         @VVV                          @VVV         \\
O(T)\otimes A^{\otimes T}   @>>>         A
\endCD $$
is commutative.

\subheading{1.3.2 Examples}
\subheading{1}
Set $\for$, $O(S)=R$ . This is called $O^{com,ass}$. Algebras over it are
commutative and associative algebras.
\subheading{2}
Set \for, $O(S)=R^{\text{all bijections:\{1,2,...,$\vert S\vert$\}}\to S}$ with
an obvious definition of composition maps. This operad (denoted $O^{ass}$)
corresponds to associative algebras.
\sbh{3}
One can define in the same manner unital operads $O^{comm,ass,1},O^{ass,1}$
and they will correspond to unital algebras.
\sbh{4}
In a similar way one defines operads $O^{Lie}$, $O^{Poisson}$, etc.
\sbh{1.3.3 Variant}
A \pm functor from an operad  $O$ to the category of
graded $\Rm$ (morphisms in this last category are
homogeneous of degree 0) will be called a graded $O$-algebra.

\subheading{1.4 Free O-algebras}

In this subsection we fix an operad $O$.
\proclaim{Lemma}
The forgetful functor ($O-algebras\to R-modules)$ admits a left adjoint.
\endproclaim
\demo{Proof}
For an $R$-module $V$ we will construct an $O$-algebra $Free(V)$:
\subheading{Construction}
$$Free(V)=\oplus_i(Free_i(V))\text{, where }Free_i(V)=(V^{\times T}
\otimes O(T))_{S^T}$$
Here $T=\{1,2,...,i\}$ and $S^T$ is the group
of permutations of the set $T$.

It is not difficult to see that $V\to Free(V)$ is the adjoint functor
we looked for.
\enddemo
\sbh{Remark}
Free $O$-algebras satisfy usual properties; e.g. if $V$ is projective
as an $R$-module then any surjection onto $Free(V)$ admits a section.
\heading 2. Modules over an algebra over an operad \endheading
\subheading{2.0}
If $A$ is an algebra over an operad $O$ we will introduce
the notion of a module over it. Our definition is motivated by a
suggestion of A.Beilinson. As it was said earlier we are dealing
only with algebras of $\Rm$
and hence modules over them will also lie in the category
of $\Rm$,
although they can be defined in a more general context.
\subheading{2.1}
Let $C$ be a \pm category and let $\Tau$ be another category. We say
that $C$ acts on $\Tau$ if $\for$ and $s\in S$ we are given a functor
$$\fss:\Csso\times{\Tau}^o\times\Tau\to\Rm$$
with $F^{\Tau}_{1,1}(\tau,\tau')=Hom(\tau,\tau')$  such that for $\sst$,
for each $\underset{t\neq t'}\to\times B_t\in C^{\times T-t'}$
and for each $\tau\in\Tau$ we
are given a natural transformation between
two functors $\Csso\times{\Tau}^o\times\Tau\to\Rm$ (composition maps): from the
functor
$${\bigotimes}_{t\in T-t'}
\fst(\underset{s\in S_t}\to\times A_s,B_t)
\otimes\fsst(\underset{s\in S_{t'}-s'}\to\times A_s,{\tau}',\tau)
\otimes F^{\Tau}_{T,t'}(\underset{t\neq t'}\to\times B_t,
\tau,{\tau}'')$$  to the functor
  $\fss(\underset{s\in S-s'}\to\times A_s,{\tau}',{\tau}'')$
\subheading{2.1.1 Examples}
\subheading{1}Any \pm category $C$ acts on itself
\subheading{2}Let $C$ be as in 1.1.1 and let $\Tau$ be a category equipped with
an action of $C$. Then when we consider $C$ as a \pm category
it will act on $\Tau$ in a natural way:
$$\fss(\underset{s\in S-s'}\to\times A_s,{\tau}_1,{\tau}_2)=
Hom(\underset{s\in S-s'}\to\otimes A_s({\tau}_1),{\tau}_2)$$

\subheading{2.1.2}
For two pairs $C_1,{\Tau}_1$ and $C_2,{\Tau}_2$ of a \pm category
and a category which it acts upon, one defines notions of \pm
functors and \pm natural transformations between \pm functors
as in 1.1.2.
\sbh{2.1.3 Variant}
One can modify the above definitions to the case of unital
\pm categories. Essentially, one needs that the unit object
$\un\in C$ "acts identically" on $\Tau$. In what follows there will
be no difference between operads and unital operads.

\subheading{2.2}
Let now $O$ be an operad and let $O$ act on a category $\Tau$.
We say that $\Tau$ is a model over $O$ if it has essentially
one object.
\subheading{2.2.1}
A model $\Tau$ can be thought of as the following linear algebra
data (analogously to 1.2 ):
\roster
\item A collection of $\Rm$ $\Tau(S,s')$ $\for$ and $s'\in S$
(they correspond to
$\fss(A^{\times(S-s')}\times\tau\times\tau)$ for $A\in O,\tau\in\Tau$)
\item A distinguished element $1\in\Tau(1,1)$
\item Composition maps $\for$ and $s'\in S$
$$\Tau(T,t')\underset{t\neq t'}\to\bigotimes O(S_t)\otimes \Tau(S_{t'},s')
\to\Tau(S,s')$$
\endroster
\subheading{2.3}
Let now $O$ be an operad and let $\Tau$ be a model over it.
A \pm functor from the pair $(O,\Tau)$ to the pair $(\Rm,\Rm)$
is called an $A$-module of type $\Tau$ for $A$ defined by the
underlying functor $O\to\Rm$. Morphisms between $A$-modules
of type $\Tau$ are defined to be natural transformations
between such functors. For any fixed $\Tau$ such modules form
an abelian category.
\subheading{2.3.1}
Again, we can describe $A$-modules of type $\Tau$ in linear
algebra terms:
An $A$-module is an $R$-module $M$ equipped with a system of maps
$\for$ and $s'\in S$
$$\Tau(S,s')\otimes A^{\otimes {S-s'}}\otimes M\to M$$ such that if
$\st$, the following diagram becomes commutative:
$$ \CD
     \trrt\otimes\arss\otimes M      @>>>  \Tau(S,s)\otimes\arss\otimes M  \\
                @VVV                                         @VVV          \\
\Tau(T,t')\otimes A^{\otimes T-t'}\otimes M       @>>>          M
\endCD $$
\subheading{2.3.2 Example}
When we put $\Tau=O$, our definition coincides with that of
[2]: $\Tau(S,s')=O(S)$,
and for $O=O^{Lie}\text{ (resp. $O^{ass,comm}$)}$ we get
usual $A$-modules (resp. Lie-algebra representations), whereas
for $O=O^{ass}$ we get $A$-bimodules. In the case of the coresponding
unital operads we get modules acted on identically by the unit.

However, by means of varying $\Tau$ the above definition alows
to get modules with an additional structure.
\sbh{2.3.3 Variant}
A \pm functor from the pair $(O,\Tau)$ to the pair (graded $\Rm$, graded
$\Rm$) will be called a graded module over the corresponding graded
algebra. We have the functor $T=shift\ of\ grading$  on the category
of graded modules over a graded algebra.

\sbh{2.4 Free modules}
We will introduce the notion of a free module over
an algebra over an operad. In particular, this will
lead to the [2] construction of the universal enveloping
algebra of an algebra over an operad.
\sbh{2.4.1}

\proclaim{Lemma}
The forgetful functor from $A$-modules of type $\Tau$
to $R$-modules admits a left adjoint.
\endproclaim
\demo{Proof}
Let $U$ be an $R$-module. We will construct an $A$-module $F(U)$,
the free $A$-module on $U$, such that the functor $U\to F(U)$
is the adjoint functor we need. Let us observe first, that
it is sufficient to construct $F(R)$ because then we can
set $F(U)=F(R){\otimes}_R U$
\sbh{Construction}
Set $F'(R)=\oplus_{i=0,1,\dots} F'_i(R)$, where
$$F'_i(R)=(\Tau(\{0,1,\dots,i\},0)\otimes A^{i})_{S^i}.$$

For any $i$ and $j$ we have maps
$$\Tau(\{0,1,\dots,i\},0)\otimes\Tau(\{0,1,\dots,j\},0)\to
\Tau(\{0,1,\dots,i+j\},0)$$
which induce on $F'(R)$ a structure of an associative
algebra with a unit. Consider now for all $i,j$ (non-commutative)
diagrams of the type:
$$ \CD
\trrt\otimes\arss\otimes F'(R) @>>> \Tau(S,s)\otimes\arss\otimes F'(R)  \\
                @VVV                              @VVV          \\
\Tau(T,t')\otimes A^{\otimes t-1}\otimes F'(R)  @>>>   F'(R)
\endCD $$
for $(S,s)'=(\{0,1,\dots,i\},0)$ and $(T,t')=(\{0,1,\dots,j\},0)$.
We define $F(R)$ as a quotient of $F'(R)$ by the ideal generated
by the images of $\phi_1-\phi_2$,
where $\phi_1$ and $\phi_2$ are two diagonal $(\searrow)$ maps
in the above diagrams. It is easy to see then, that
$F(R)$ constructed in this way satisfies the properties we need.
\enddemo
\sbh{2.4.2}
Put now $U=R$ and let us donote by $P_A$ the corresponding $A-$module $F(R)$:
$Hom_A(P_A,M)\simeq M$ as an $R$-module, for any $A$-module $M$.
Then $P_A$ has a natural structure of an associative algebra
with a unit, since $P_A\simeq End_A(P_A)$, and the category of $A-$
modules of type $\Tau$ is naturally equivalent to the category
of right $P_A$-modules. For $\Tau=O$ the algebra $P_A$ is
the universal enveloping algebra of $A$ in the terminology of [2].
\sbh{2.4.3}
Let now $B$ be another $O$-algebra and let us have
a homomorphism from $B$ to $A$. Then the obvious restriction
functor from the category of $A$-modules to the category
of $B$-modules admits a left adjoint, called the iduction
functor. Its existence is obvious from the equivalence
of categories of 2.4.2 and the fact that we have an associative
algebra homomorphism $P_B\to P_A$.

\subheading{2.5 Derivations}
{}From now on we will consider modules over an algebra over
an operad with $\Tau=O$ and we will call them just $A$-modules.
\sbh{2.5.1}
Let $A$ be an olgebra over an operad $O$ and let $M$ be an
$A$-module. \newline
An $R$-linear map $\phi:A\to M$ is said to be a derivation
(from $A$ to $M$) if \newline
$\for$  the diagram
$$ \CD
O(S)\otimes\as   @>\sre({\phi}^{\otimes S-s}\otimes id)>>  \sre(O(S,s)\otimes
A^{\otimes S-s}\otimes M) \\
  @VVV                                                             @VVV
\\
   A                               @>\phi>>                          M
\endCD $$
is commutative.

The set of all derivations from $A$ to $M$ will be denoted by
$\Omega (A,M)$

\subheading{2.5.2}
Suppose that $A$ is a free algebra $A=Free(V)$. Then $\Omega(A,M)=
Hom_R(V,M)$ for any $A$-module $M$.

\heading 3. Sheaves and cohomology \endheading

\subheading{3.0}
Starting with an $O$-algebra $A$, we will construct a site $C(A)$.
This definition appeared first in [1] where Quillen proved that
cohomologies of certain sheaves on this site provide correct
cohomology theories for $A$-modules (in the case of commutative,
associative and Lie algebras).
\subheading{3.1}
Let us consider the category $C(A)$ consisting
of $O$-algebras $B$ with a
homomorphism $B\to A$. This category possesses a fibered product
and we introduce a Grothendieck topology on it by declaring
epimorphisms to be the covering maps.

Thus we can consider sheaves on $C(A)$ and their cohomologies.
\subheading{3.2 Sheaves ${\Im}_M$}

Let an $M$ be an $A$-module. Then it is also a module over each
algebra $B\in C(A)$. We define a sheaf $\fm$  on $C(A)$ by setting
$$\Gamma(B,\fm)=\Omega(B,M)$$
Sheaf axioms are easily verified.
\sbh{3.2.1}
The following remark is due to essentially to [1]:
Each $\fm$ is representable by $A_M\in C(A)$, equal to
$A\oplus M$  with the natural algebra structure on it,
in other words $$\Gamma(B,\fm)=\text{algebra homomorphisms }B\to A\oplus M.$$
It is also not difficult to see that the functor $M\to\fm$ is fully
faithful.
\sbh{3.3}
Let us mention several properties of the
category $C(A)$.
\sbh{3.3.1}
If $E\to D$ is a covering in $C(A)$, then the functor \newline
$C(A)_D\text{ (objects of $C(A)$ over $D$)}\to\text
{descent data on $E$ with respect to $D$}$ \newline
is an equivalence of categories.
\sbh{3.3.2}
If $V$ is a projective $R$-module with an $R$-module map to $A$, then \newline
$Free(V)\in C(A)$ and the functor $\Im\to\Gamma(Free(V),\Im)$
is exact. This follows e.g. from Remark 1.4.
\sbh{3.3.3}
\proclaim{Proposition}
The functor $\Im$: ($A$-modules $\to$ sheaves) admits right and left
adjoint functors, $\RR$ and $\LL$ respectively. We have $\RR\circ\Im\simeq
\LL\circ\Im\simeq Id_{A-mod}$.
\endproclaim
\demo{Proof}
For each $X\in C(A)$ consider the sheaf $Const_X$ defined by
$$\Gamma(Y,Const_X)=\underset Hom(Y,X)\to\oplus R.$$
The sheaf $Const_X$ is defined uniquely by the following property:
$Hom(Const_X,S)=\Gamma(X,S)$ functorially in $S-$a sheaf over $C(A)$.
Let $B=Free(V)$. Then $$Hom(Const_{Free(V)},\fm)=\Gamma
(Free(V),\fm)=\Omega(Free(V),M)=Hom_R(V,M).$$
This fact together with 3.3.2 imply that the functor $\Im$ is exact.

In order to construct the functor $\LL$, it is sufficiant
to define the values of $\LL$ on sheaves of the form
$Const_B$ for free algebras $B=Free(V)$, because any sheaf over $C(A)$
is a quotient of a direct sum of such sheaves.
However, the above calculation shows that for these sheaves we can put
$\LL(Const_{Free(V)})=F(V)$ in the notation of 2.4.

To construct the functor $\RR$ we put for a sheaf $\FF$,
$\RR(\FF)=Hom(\Im(P_A),\FF)$ with the obvious structure
of a right module over $P_A$. Then we use 2.4.2.

The fact that $\RR\circ\Im\simeq\LL\circ\Im\simeq Id_{A-mod}$
follows from the full faithfulness of the functor $\Im$.
\enddemo
\sbh{3.3.4}
As always, the functor $\LL$ is right exact and the functor
$\RR$ is left exact and we can consider their left (resp. right)
derived functors $L^{\cdot}\LL$ (resp. $R^{\cdot}\RR$).
(The functor $\LL$ can be derived e.g. because
the sheaves $Const_{Free(V)}$ for projective $V$ are projective (!)
in the category of sheaves of $\Rm$.)
\sbh{3.4}
Let us apply the functor $L^{\cdot}\LL$ to the sheaf $R_A$.
We obtain an object $T_A$ in $D(A)$
(the derived category of $A$-modules). $T_A$ is called the
cotangent complex of $A$.
$$RHom(T_A,M)\simeq R\Gamma(A,{\Im}_M)$$
\sbh{3.4.1}
The fact that $\LL$ is right exact implies that
$$H^i(T_A)=\cases
0, & i>0 \cr
I_A, & i=0
\endcases$$
where $I_A$ is the $A$-module representing the functor
$M\to\Omega(A,M)$.
\sbh{3.4.2 Examples}
\sbh{1} If $O$ is the Lie operad $O^{Lie}$, $I_A$ canonically identifies
with the augmentation ideal of the universal enveloping
algebra.
\sbh{2} If $O=O^{ass,1}$, $I_A\simeq I=ker(A\otimes A\to A)$.
\sbh{3} If $O=O^{com,ass,1}$, $I_A\simeq I/I^2$, with
$I$ is as above.

\heading 4. Deformations \endheading
\sbh{4.0}
Results of this section are partially contained in [1],[3] and [4].
We decided to present them, since the formalism developed in the
preceeding sections seems to be a covenient tool for passage
from deformations to cohomology. By definition, we put
$H^i(A,M)=RHom^i(T_A,M)$.
\sbh{4.1}
Let $A$ be an $O$-algebra (over $R$).
An $i$-th level deformation of $A$ is an $O$ algebra
$A_i$ over $\rti$ (in the sense of 1.2.4) with an isomorphism
$\phi :A_i/t\cdot A_i\simeq A$ and such that
$$Tor_1^{\rti}(A_i,R)=0$$
In other words, we need that
$$ker(t:A_i\to A_i)=im(t^i:A_i\to A_i)\text{ identifies under
a natural map with } A$$
\sbh{4.1.1}
The category of $i$-th level deformations (morphisms are compatible
with $\phi$'s) is a groupoid denoted by $Deform^i(A)$. For each $i$ there
are functors from $Deform^{i+1}(A)$ to $Deform^i(A)$ (taking modulo $t^{i+1}$).
If $A_i\in Deform^i(A)$, the fiber $\Cai$ of $Deform^{i+1}(A)$ over $A_i$
will be called the category of prolongations of $A_i$.
\sbh{4.2}
\proclaim{Theorem}
\roster
\item
The category of $1$-st level deformations is equivalent
to the category $T(\fa)$
of $\fa$-torsors. In particular, $\pi_o(C_1)\simeq
H^1(A,A)$
\item
If $A_{i+1}\in\Cai$, $Aut(A_{i+1})$ as of an object of this
category is canonically isomorphic to $\Omega(A,A)$
\item
To each $A_i\in Deform^i(A)$ we can associate a gerbe $G_{A_i}$ over
$C(A)$ bounded by $\fa$ in such a way that $G_{A_i}$ is canonically
equivalent to $\Cai$.

In particular, this means that to each $A_i$ we can assign
an element in $H^2(A,A)$ that vanishes if and only if
$\Cai$ is nonempty. And if $\Cai$ is nonempty, its $\pi_o$
is a torsor over $H^1(A,A)$.
\endroster
\endproclaim
\demo{Proof}
\roster
\item
The functor $T:C(A)\to T(\fa)$ is given by:
$$\Gamma(B,T(A_1))=O-\text{algebras homomorhisms over R:}B\to A_1$$
Using 3.3.1 it is easy to show that it is an equivalence
of categories.
\item
This is a direct verification.
\item
We define the gerbe as follows:

$G_{A_i}(B)$ is the groupoid of $R[t]/t^{i+2}\cdot R[t]-O$ algebras $B_{i+2}$
with an isomorphism  $B_{i+1}/t^{i+1}\cdot B_{i+1}\simeq
B\underset{A}\to\times A_i$ such that \newline
$ker(t^{i+1}:B_{i+1}\to B_{i+1})=im(t:B_{i+1}\to B_{i+1})$ and
identifies under a natural morphism with $A_i$.

Functors $G_{A_i}(B)\to G_{A_i}(D)$ for maps $D\to B$ are given
by taking fibre products.

It is easy to check that $G_{A_i}$ is indeed a gerbe bounded by
$\fa$ over $C(A)$ and that its fiber over $A$ is equivalent to $\Cai$.
\endroster
\enddemo
\sbh{4.2.1} To summarize, we have shown that the cohomology groups
$H^i(A,A)$ "control" the deformation theory of $A$.
\sbh{4.3}
For the remainder of this paper we restrict ourselves
to the case $O=O^{ass,1}$.

\sbh{4.3.1}
When $A$ is flat as an $R$-module, we have a theorem of Quillen [1]:
\proclaim{Theorem}
$H^i(T_A)=0$ for $i\neq 0$.
\endproclaim
By the cohomology long exact sequence of the triple
$$0\to I_A\to A\otimes A\to A\to 0$$ we get that in this case
$H^i(A,M)=\Ei(A,M)$ for any $A$-bimodule $M$ and $i\ge 1$.
\sbh{Variants}
\sbh{4.3.2}
If $A$ is an augmented algebra we can look for its
deformations in the class of augmented algebras.
Then Theorem 4.2 remains valid after replacing $H^{\cdot}(A,A)$ by
$H^{\cdot}(A,A_+)$, where $A_+$ denotes the augmentation ideal of $A$.

\sbh{4.3.3}
Suppose now that $A$ is a graded algebra. Then we will consider
the site $C(A)$ that corresponds to graded algebras over $A$.
If now $M$ is a graded $A$-bimodule, we introduce graded cohomology groups as
$(H^i(A,M))_j=H^i(A,T^j(M))$, with $T$ being the translation functor of 2.3.3.

A graded deformation of $A$ of $i$-th level is an algebra
$A_i$ as above endowed with a grading such that $deg(t)=1$.
Then the Theorem 4.2 is restated in the following way:
\roster
\item
Isomorphism classes of first level deformations are in 1-1
correspondence with the elements of $(H^1(A,A))_{-1}$
\item The automorphism group of a prolongation of an $i$-th
level deformation is naturally isomorphic to $\Omega(A,A)_{-i-1}$
\item The obstruction to the existence of a
prolongation of a given $i$-th level deformation
lies in $(H^2(A,A))_{-i-1}$
\item The set of isomorphism classes of prolongations of
a given $i$-th level deformation is a torsor over $(H^1(A,A))_{-i-1}$
\endroster

In the graded case Theorem 4.3.1 states that
$$(H^i(A,M))_j=(\Ei(A,M))_j$$ for any graded $A$-bimodule
$M$.
\sbh{4.3.3'}
Suppose now that we have a family $A_i$ of graded deformations
which are prolongations of one another. In this case we can form
an algebra
$$A_t=\text{elements of finite degree in }\underset{\longleftarrow}\to {lim}
(A_i)$$

This is an algebra over $R[t]$.
Consider its fiber at $t=1:A_1=A_t/(t-1)\cdot A_t$. This algebra will carry
a natural filtration and the associated graded algebra $gr(A_1)$
will be canonically isomorphic to $A$.
\sbh{4.3.4}
One, of course, can consider a combination of graded and augmented
situations. Statement of Theorem 4.2 will change correspondingly.

\heading 5. The Poincar\'e-Birkhoff-Witt theorem \endheading
\sbh{5.0}
In this section
we will give an application of the theory presented above. We will prove
the \pbw  for Lie algebras over a ring $R$ which are flat
as $\Rm$.
\sbh{5.1}
Let us recall main definitions.
\sbh{5.1.1}
Let $\g$ be an $R$-module. Its $i$-th exterior power $\Lamda ^i(\g)$
is defined to be the subspace of $\g ^{\otimes i}$ spanned
by tensors of the form $Alt(g_1,g_2,\dots,g_i)$, where $Alt$
means alternating sum.
\sbh{5.1.2}
$S^{\cdot}(\g)$ will denote the symmetric algebra of $\g$ which
is by definition the quotient of the tensor algebra $T^{\cdot}(\g)$
by the ideal generated by $\Lamda ^2(\g)$. It is a graded algebra
with graded components denoted by $S^i(\g)$.
\sbh{5.1.3}
A Lie algebra structure on $\g$ is a map $[,]:\Lamda ^2(\g)\to\g$
such that the map $[,]\circ (id\otimes [,]-[,]\otimes id):\Lamda ^3
(\g)\to\g$ vanishes.
\sbh{5.1.4}
For a Lie algebra $\g$ its universal enveloping algebra $U(\g)$
is the quotient of the tensor algebra $T^{\cdot}(\g)$ by the
ideal generated by $\omega -[,](\omega)$ for $\omega\in\Lamda ^2(\g)$.
$U(\g)$ caries a natural filtration and there is a canonical
surjection $S^{\cdot}(\g)\to gr(U(\g))$.
\sbh{5.2 The \pbw }
\proclaim{Theorem}
Let $\g$ be a Lie algebra over $R$ which is flat as an $R$-module.
Then the canonical epimorphism $S(\g)\to gr(U(\g)$ is an
isomorphism.
\endproclaim
\sbh{5.2.1}
We will prove this theorem by constructing a graded deformation
as in 4.3.4 of $S(\g)$ in the class of augmented algebras. This idea
is borrowed from [5], where the \pbw is proved
for Koszul algebras also by means of deformation theory.

By 4.3.2 and 4.3.4 we know that in this case the
deformations of $A=S(\g)$ are controled by
$Ext_{A\otimes A}(A,A_{+})$'s.
\sbh{5.2.2}
\proclaim{Proposition}
Let $\g$ be a flat $R$-module.
Then
\roster
\item
$$(Ext_{A\otimes A}^2)_{-j}=\cases
0, & j>1 \cr
Hom_R(\Lamda^2(\g),\g) if & j=1
\endcases$$
\item
$$(Ext_{A\otimes A}^3)_{-j}=\cases
0, & j>2 \cr
Hom_R(\Lamda^3(\g),\g) if & j=2
\endcases $$

\endroster
\endproclaim

This proposition easily follows from the fact
that for flat $\g$-modules the Koszul complex
$$\dots\to S(\g)\otimes\Lamda^i(\g)\otimes S(\g)\to\dots\to
S(\g)\otimes\Lamda^2(\g)\otimes S(\g)\to S(\g)\otimes\g\otimes S(\g)
\to S(\g)$$
is exact.
\sbh{5.3}
\demo{Proof of the theorem}
We put $A=S(\g)$. This is a graded augmented algebra
and by a deformation we will mean a graded deformation in
the class of augmented algebras.
\sbh{Step 1}
Let us build a first level deformation of $A$ that corresponds to
the cohomology class of our $[,]:\Lamda ^2(\g)\to\g$ via 4.3.4
and 5.2.2. It is canonical up to a unique
isomorphism since $\Omega(A,(A_{+})_{-1}=0$.

Untwisting of the passage ( first level deformations $\to$
cohomology) shows that there exists a canonical $R$-linear
isomorphism $\phi:\g\to ((A_1)_{+})_1$ such that
$$\phi(x)\cdot\phi(y)-\phi(y)\cdot\phi(x)=t\phi([x,y])$$
\sbh{Step 2}
A direct verification shows that the obstruction to
the existense of a prolongation of this deformation as
a map $\Lamda ^3(\g)\to\g$ is given by the expression
$[,]\circ (id\otimes [,]-[,]\otimes id)$ which in turn vanishes by the
Jacobi identity.
So a prolongation exists and is unique up to a unique
isomorphism by 4.3.4. and 5.2.2.
\sbh{Step 3}
Again by 4.3.4  and 5.2.2 for any $i\geq 2$ an $i$-th level deformation
of $A$ can be prolonged in a unique up to a unique isomorphism
way and so we find ourselves in the situation of 4.3.3'.
\sbh{Step 4}
\proclaim{Claim}
There exists a canonical $R$-linear map $\phi':\g\to ((A_t)_{+})_1$ such that
$$\phi'(x)\cdot\phi'(y)-\phi'(y)\cdot\phi'(x)=t\phi'[x,y]$$
\endproclaim
Ideed, since $deg(t)=1$, there is a unique way to lift $\phi$
of Step 1 to a map $\g\to ((A_t)_{+})_1$.

Then, by the definition of $U(\g)$ there exists a map
$\phi':U(\g)\to A_1$ that prolongs the map $\phi'$ above.
\sbh{Step 5}
Let us consider the associated graded map
$gr(\phi'):gr(U(\g)\to gr(A_1)$ and let us also consider
the composition
$$S(\g)\underset{\text{by 5.1.4}}\to\to gr(U(\g))\underset
{\phi'}\to\to gr(A_1)\underset{\text{by 4.3.3'}}\to\simeq S(\g)$$

This composition is easily seen to be the identity map. This implies that
\roster
\item
$\phi'$ is an isomorphism between $U(\g)$ and $A_1$
\item
The \pbw.

\endroster
\enddemo

\enddocument